\documentclass{article}
\newcommand{\CG}{\mbox{$\mathcal{G}$}} 
\newcommand{\CP}{\mbox{$\mathcal{P}$}}
\newcommand{\CQ}{\mbox{$\mathcal{Q}$}}
\newtheorem{theorem}{Theorem}
\newtheorem{definition}{Definition}
\usepackage{graphicx}
\usepackage{amssymb}

\begin{document}

\title{Semisymmetric Graphs from Polytopes} 

\author{Barry Monson\thanks{Supported by the NSERC of Canada, grant 4818}
\\University of New Brunswick
\\Fredericton, New Brunswick, Canada E3B 5A3 (bmonson@unb.ca)
\and Toma\v{z} Pisanski\thanks{Supported 
by Ministry of Higher Education, Science and Technolgy of
Slovenia grants P1-0294,J1-6062,L1-7230.}
\\ University of Ljubljana, FMF
\\ Jadranska 19, Ljubljana 1111, Slovenia; 
\\ and  Faculty of Education, University of Primorska, 
\\ Cankarjeva 5, Koper 6000,  Slovenia
(tomaz.pisanski@fmf.uni-lj.si)
\and Egon Schulte\thanks{Supported 
by NSA-grant
H98230-05-1-0027}
\\ Northeastern University
\\ Boston MA 02115, USA (schulte@neu.edu)
\and Asia Ivi\'{c} Weiss\thanks{Supported by the NSERC of Canada, grant 8857}
\\York University
\\Toronto, Ontario, Canada M3J 1P3 (weiss@yorku.ca)}
\date{}
\maketitle

\vspace{5mm}
\begin{abstract}
\noindent
Every finite, self-dual, regular (or chiral) $4$-polytope
of type $\{3,q,3\}$ has a trivalent $3$-transitive (or $2$-transitive)
medial layer graph. Here, by dropping self-duality, we obtain
a construction for semisymmetric trivalent graphs 
(which are edge-  but not vertex-transitive). 
In particular, the Gray graph arises
as the medial layer graph of a certain universal locally toroidal
regular $4$-polytope. 



\end{abstract}
\vspace{5mm}
%
%



\section{Introduction}

The theory of symmetric trivalent graphs and   the theory 
of regular polytopes
are each abundant sources of beautiful mathematical ideas.
In \cite{mowe}, two of the authors established some general and unexpected
connections between the two subjects, building upon a rich variety
of examples appearing in the literature (see \cite{fost}, \cite{cond2}, 
\cite{cox1},
\cite{coxweiss}, \cite{triv} and \cite{infg}). 
Here we develop these connections a little further, 
with specific focus on semisymmetric graphs. In particular, we reexamine
the \textit{Gray graph},  described in \cite{bou1,bou2} and \cite{marpar,pisaram},
and here appearing as the medial layer graph of an 
abstract regular $4$-polytope.

We begin  with some basic ideas concerning symmetric
graphs \cite[ch. 18-19]{biggs}.  Although some  of the following results
generalize to graphs of higher valency,  for brevity we shall 
{\em assume outright that $\mathcal{G}$ is a simple, finite, connected trivalent 
graph}\, (so that each vertex has valency $3$).

By a {\em $t$-arc} in $\mathcal{G}$ we mean a list of
vertices $[v] = [v_0,v_1,\ldots,v_t]$ such that $\{v_{i-1},v_i\}$ is
an edge for $1 \leq i \leq t$, but no $v_{i-1} = v_{i+1}$.
Tutte has shown that there exists a maximal
value of $t$ such that the automorphism group ${\rm Aut}(\mathcal{G})$ 
is transitive on $t$-arcs. 
We say that
$\mathcal{G}$ is
\textit{$t$-transitive} if ${\rm Aut}(\mathcal{G})$ is transitive on $t$-arcs, 
but not
on $(t+1)$-arcs in $\mathcal{G}$, for some $t\geq 1$.
Tutte also proved the remarkable result that   
$t \leq 5$ (\cite[Th.18.6]{biggs}). Any such
arc-transitive graph is said to  be \textit{symmetric}.

Each fixed $t$-arc $[v]$ in  a $t$-transitive graph $\mathcal{G}$ 
has {\em stabilizer sequence}   

\begin{displaymath}
\mbox{Aut}\;\mathcal{G} \supset B_t \supset B_{t-1} \supset \ldots \supset B_1 
\supset B_0,
\end{displaymath}
where the subgroup $B_j$ is the pointwise stabilizer of $\{v_0,
\ldots, v_{t-j}\}$.
Since   ${\rm Aut}(\mathcal{G})$ is transitive on $r$-arcs, for $r \leq t$, 
the subgroup
$B_j$  is conjugate to that obtained from any other $t$-arc.  
In particular,
$B_t$ is the  vertex stabilizer, whereas $B_0$
is the pointwise stabilizer of the whole arc.   In fact, 
$B_0 = \{ \epsilon\}$ is trivial (\cite[Prop. 18.1]{biggs}), so that 
${\rm Aut}{(\mathcal{G})}$ acts sharply transitively on $t$-arcs.

Each $t$-arc $[v]$ has two \textit{successors}, $t$-arcs of the form
$[v^{(k)}] := [v_1, \ldots, v_t, y_k]$, where $v_{t-1}, y_1, y_2$ are
the vertices adjacent to $v_t$. The \textit{shunt} $\tau_k$ is the 
(unique) automorphism of $\mathcal{G}$ such that
$[v] \tau_k = [v^{(k)}]$. Also let $\alpha$ be the unique
automorphism which reverses the basic
$t$-arc $[v]$. Then $\alpha$ has period $2$ and $\alpha \tau_1 \alpha$
equals either  $\tau_1^{-1}$ or $\tau_2^{-1}$. We shall say that
$\mathcal{G}$ is of \textit{type} $t^+$ or $t^-$, respectively.
We can now assemble several beautiful
results concerning  ${\rm Aut}{(\mathcal{G})}$ (see  \cite[ch. 18]{biggs}).  

\begin{theorem}\label{autgraph}
 Suppose $\mathcal{G}$ is a finite connected
$t$-transitive trivalent graph, with $1 \leq t$, and
suppose $\mathcal{G}$ has $N$  vertices.  Then 

{\em (a)} For $0 \leq j \leq t-1$ we have $|B_j| = 2^j$. Also,
$|B_t| = 3 \cdot 2^{t-1}$ and 
$$|\mbox{Aut} (\mathcal{G})| = 3 \cdot N \cdot 2^{t-1}\;.$$

{\em (b)} The stabilizers $B_j$ are determined up to isomorphism by $t$\,\rm{:}

$$ \begin{array}{c|c|c|c|c|c}
t & B_1 & B_2 & B_3 & B_4 & B_5\\ \hline
1 & \mathbb{Z}_3 &  & & &\\ \hline
2 & \mathbb{Z}_2 & \mathbb{S}_3 &  &  & \\ \hline
3 & \mathbb{Z}_2 & ( \mathbb{Z}_2 )^2 & \mathbb{D}_{12} &   &\\ \hline
4 & \mathbb{Z}_2 & ( \mathbb{Z}_2)^2 & \mathbb{D}_8 & \mathbb{S}_4 &   \\
\hline
5 & \mathbb{Z}_2 & ( \mathbb{Z}_2)^2 & (\mathbb{Z}_2)^3 & \mathbb{D}_8
\times \mathbb{Z}_2 & \mathbb{S}_4 \times \mathbb{Z}_2 \\ \hline
\end{array} $$

{\rm (c)} $\mathcal{G}$ is one of $7$ types: 
$1^-, 2^+, 2^-, 3^+, 4^+, 4^-$ or\, $5^+$.
\end{theorem} 
\noindent
(Here $\mathbb{Z}_k$ is the cyclic group of order $k$, $\mathbb{D}_{2k}$
is the dihedral group of order $2k$, $\mathbb{S}_k$ is the symmetric
group of degree $k$.)

 Useful lists of symmetric trivalent graphs appear in
\cite{fost} and \cite{cond2}. We refer to \cite{mowe} for a  
description of several interesting examples.

\medskip

We   now  briefly describe some key properties of
abstract regular and chiral polytopes, referring again to
\cite{mowe} for a short discussion, and to \cite{arp,schwei,scwe} for details.
An {\em (abstract) $n$-polytope} $\mathcal{P}$ is a partially ordered set
with a strictly monotone rank function having range $\{-
1,0,\ldots,n\}$. An element $F \in \mathcal{P}$ with ${\rm rank}(F)=j$ is
called a $j$-{\it face}; typically $F_j$ will indicate a $j$-face; and 
$\mathcal{P}$ has a unique least face
$F_{-1}$ and   unique greatest face $F_n$. 
Each maximal
chain or {\it flag} in $\mathcal{P}$ must contain $n+2$ faces. Next, 
$\mathcal{P}$ must satisfy a 
\textit{homogeneity property}\,: whenever $F < G$ with ${\rm rank}(F)=j-1$ and
${\rm rank}(G)=j+1$, there are exactly two $j$-faces $H$ with
$F<H<G$, just as happens for convex $n$-polytopes. 
It follows that for $0 \leq j \leq n-1$ and any flag $\Phi$, 
there exists a unique \textit{adjacent} 
flag $\Phi^j$, differing from $\Phi$ in just the
rank $j$ face. With this notion of adjacency the flags of $\mathcal{P}$
form a \textit{flag graph}\, (not to be confused with the medial layer 
graphs appearing below).
The final defining property of $\mathcal{P}$ is that it  should
be \textit{strongly flag--connected}. This means that the flag graph for each
section is connected. 
Whenever $F \leq G$ are faces of ranks $j \leq k$
in  $\mathcal{P}$,  the \textit{section}
$ G/F := \{ H \in \mathcal{P}\, | \, F \leq H \leq G \}$
is thus in its own right a  ($k-j-1$)-polytope.

Since our main concern is with 4-polytopes, we
now tailor our discussion to that case.  A (rank 4) polytope $\mathcal{P}$
is  {\em equivelar} of type $\{p_1,p_2,p_3\}$ if, for $j =
1,2,3$, whenever $F$ and $G$ are incident faces of $\mathcal{P}$ with
$\mbox{rank} (F) = j-2$ and $\mbox{rank} (G) = j+1$, then the rank 2
section $G/F$ has the structure of a $p_j$-gon (independent of choice
of $F < G$).  Thus, each 2-face (polygon) of $\mathcal{P}$ is isomorphic to
a $p_1$-gon, and there are $p_3$ of these arranged around each 1-face
(edge) of $\mathcal{P}$; and in every 3-face (facet) of $\mathcal{P}$, each
0-face is surrounded by an alternating cycle of $p_2$ edges and $p_2$
polygons.

The automorphism group $\mbox{Aut}(\mathcal{P})$ consists of all
order preserving bijections on $\mathcal{P}$. If $\mathcal{P}$
also admits a duality (order reversing bijection), then $\mathcal{P}$
is said to be \textit{self-dual}; clearly $\mbox{Aut}(\mathcal{P})$
then has index 2 in the group $D(\mathcal{P})$ of all automorphisms and
dualities.  (Note that $D(\mathcal{P}) =
\mbox{Aut}(\mathcal{P})$ when $\mathcal{P}$ is not self-dual.) If $\mathcal{P}$
is self-dual and equivelar, then it has type $\{p_1 , p_2 ,p_1\}$.

\begin{definition}\label{meddef} Let $\mathcal{P}$ be a 4-polytope.  
The associated
{\em medial layer graph} $\mathcal{G}(\mathcal{P})$, or briefly $\mathcal{G}$, is the
simple graph whose vertex set is comprised of all 1-faces and 2-faces
in $\mathcal{P}$, two such taken to be adjacent when incident in $\mathcal{P}$.
\end{definition}

\noindent\textbf{Remarks}:
Any medial layer graph $\mathcal{G}$ is easily seen to be 
bipartite and connected. Note that the more desirable phrase
`medial graph' already has a somewhat different meaning in the 
literature on topological graph theory.

\medskip

To further focus our 
investigations, we \textit{henceforth assume that $\mathcal{P}$ is equivelar
of type $\{3,q,3\}$}, where the integer $q\geq 2$. Thus
$\mathcal{G}$ is trivalent, with vertices of two types occuring  
alternately along cycles of length $2q$.  
We  say that a $t$-arc in $\mathcal{G}$
is of {\em type $1$} (resp. {\em type $2$}\,) if its 
initial vertex is a $1$-face (resp. $2$-face) of $\mathcal{P}$.
The fact that certain polygonal sections
of $\mathcal{P}$ are triangular immediately implies that the action 
of $D(\mathcal{P})$ on $\mathcal{G}$ is faithful, so that we 
may regard $D(\mathcal{P})$, or
$\mbox{Aut}(\mathcal{P})$, as  a subgroup of $\mbox{Aut}(\mathcal{G})$
(see \cite[\S\,2]{mowe}).

\begin{figure}[hbt]
\begin{center}
\includegraphics[width=105mm]{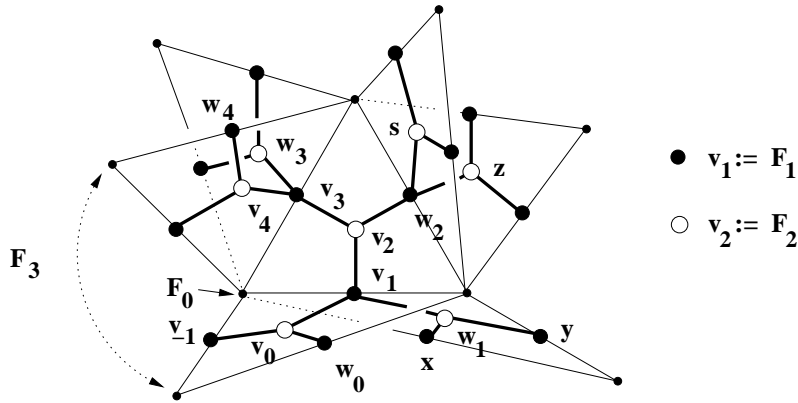}
\end{center}
\caption{A fragment of a polytope of type $\{3,q,3\}$, with $q=6$.}
\label{frag}
\end{figure}

In Figure \ref{frag} we show a fragment of a polytope $\mathcal{P}$ of
type $\{3,6,3\}$.  The vertices of $\mathcal{G}$ are here represented as
black and white discs, and the edges of
$\mathcal{G}$ are indicated by heavy lines. 

\medskip 

Since we shall soon assume that $\mathcal{P}$ is quite symmetric,
it is useful now to fix 
a \textit{base flag} 
$$\Phi = \{F_{-1},F_0, F_1, F_2, F_{3}, F_4\} $$ 
in $\mathcal{P}$. Given this,  
it is convenient to define
$v_1 := F_1, v_2 := F_2$,
and in general let $v_0 = v_{2q} , v_1, v_2, \ldots ,v_{2q -1} = v_{-1}$
denote alternate edges and polygons in the rank $2$ section
$F_3 / F_0$ of $\mathcal{P}$. Thus each $v_j $ is adjacent in 
$\mathcal{G}$ to $v_{j \pm 1}$,
taking subscripts mod $2q$. We also let $w_j$ be the third vertex adjacent to 
$v_j$ in $\mathcal{G}$, as indicated in Figure~\ref{frag}.

We turn now to two significant classes of highly symmetric
polytopes. First we recall that  
$\mathcal{P}$ is  
\textit{regular} when  $\rm{Aut}(\mathcal{P})$ acts 
transitively on the flags
of $\mathcal{P}$. Assuming still that $n=4$, we observe
that for $0 \leq j \leq 3$, there exists a 
(unique) automorphism $\rho_j$ mapping the base flag $\Phi$ 
to the adjacent 
flag $\Phi^j$. Then
$\rm{Aut}(\mathcal{P})$ is generated by the involutions  
$\rho_0,\rho_1,\rho_2,\rho_3$, which 
satisfy  at least the relations 

\begin{equation}\label{regrels}
\begin{array}{c}
\rho_0^2 = \rho_1^2 = \rho_2^2 = \rho_3^2 = (\rho_0 \rho_2)^2 =
(\rho_0 \rho_3)^2 = (\rho_1 \rho_3)^2 = \epsilon\\
  \\
(\rho_0 \rho_1)^{p_1} =
(\rho_1 \rho_2)^{p_2} = (\rho_2 \rho_3)^{p_3} = \epsilon\;\;,
\end{array}
\end{equation}
with $2 \leq p_1 , p_2 , p_3  \leq \infty$. Indeed, $\mathcal{P}$ is  
equivelar
of type $\{p_1 , p_2 , p_3\}$. (As before we will actually have
$p_1 = p_3= 3$ and $p_2 = q$ in our applications.)

Furthermore, an \textit{intersection
condition} on standard subgroups holds: 
\begin{equation}
\label{reginter}
\langle \rho_i\,|\, i \in I \rangle \cap \langle \rho_i\,|\,i \in J \rangle =
\langle \rho_i\,|\, i \in I \cap J \rangle
\end{equation}
for all $I,J \subseteq \{0,1,2,3\}$. 
In short, ${\rm Aut}(\mathcal{P})$
is a very particular
quotient of a Coxeter group with string diagram. 
 
Conversely,  suppose that $\Gamma = \langle \rho_0,\ldots,\rho_{3}
\rangle$ is a \textit{string C-group}, namely any group generated by 
specified involutions 
satisfying (\ref{regrels}) and (\ref{reginter}).  
Then one may construct a
regular $4$-polytope
$\mathcal{P} = \mathcal{P}(\Gamma)$,  of type $\{ p_1,p_2,p_3\}$,
with ${\rm Aut}(\mathcal{P}) = \Gamma$.
We refer to \cite[Def. 2]{mowe} or \cite[Thms. 2E11 and 2E12]{arp} for 
details of the construction. Note also that 
$\mathcal{P}$ is self-dual if and only if 
${\rm Aut}(\mathcal{P})$ admits an involutory 
group automorphism $\delta$ such that
$\delta \rho_j \delta = \rho_{3-j}$ for $j = 0,1,2,3$.
Such a polytope $\mathcal{P}$ admits a \textit{polarity} 
(i.e. involutory duality)
which reverses the basic flag $\Phi$.  Thus 
$D(\mathcal{P}) \simeq {\rm Aut}(\mathcal{P})  \rtimes \mathbb{Z}_2$  
(see \cite[2B17 and 2E12]{arp}).

The upshot of Theorem 2 in \cite{mowe}
is that $\mathcal{G}(\mathcal{P})$ is $3$-transitive 
when $\mathcal{P}$ is finite, regular
and self-dual of type $\{3,q,3\}$.

\medskip

For any regular polytope $\mathcal{P}$, the \textit{rotations}
$\sigma_j:=\rho_{j-1}\rho_j$  generate a  subgroup
${\rm Aut}(\mathcal{P})^+$ having index 1 or 2 in ${\rm Aut}(\mathcal{P})$. 
In the
latter case, $\mathcal{P}$ is said to be {\it directly regular}, and certain
properties of the $\sigma_j$ lead, in a natural way, to 
a parallel theory of  {\it chiral}
polytopes (see \cite{schwei,scwe} for details).

A polytope $\mathcal{P}$ of rank $n \geq 3$ is said to be 
\textit{chiral} if it is not regular,  but there do
exist automorphisms $\sigma_1,\ldots,\sigma_{n-1}$ such
that $\sigma_j$ fixes all faces in $\Phi \backslash \{ F_{j-1},F_j\}$
and cyclically permutes consecutive $j$-faces of $\mathcal{P}$ in the
rank 2 section $F_{j+1}/F_{j-2}$ of $\mathcal{P}$. 
The automorphism group of $\mathcal{P}$ now has two flag
orbits, with adjacent flags always in different orbits.
Again taking $n = 4$, it is even possible in the chiral case to choose
automorphisms $\sigma_1, \sigma_2, \sigma_3$ which generate 
${\rm Aut}(\mathcal{P})$ and satisfy at least the relations
\begin{equation}
\label{chirrels}
\begin{array}{c}
\sigma_1^{p_1} = \sigma_2^{p_2} = \sigma_3^{p_3} = \epsilon\\
(\sigma_1 \sigma_2)^{2} = (\sigma_2 \sigma_3)^{2} = 
(\sigma_1 \sigma_2 \sigma_3)^{2} = \epsilon,
\end{array}
\end{equation}
for some $2 \leq p_1 , p_2 , p_3 \leq \infty$. Once more $\mathcal{P}$ is 
equivelar of
type $\{ p_1, p_2, p_3 \}$. Here too the specified 
generators satisfy a revised intersection condition:

\begin{equation}
\label{chirinter}
\begin{array}{rcl}
	  \langle \sigma_1\rangle \cap \langle \sigma_2 \rangle
	&= \: \{\epsilon\} \: =& \langle \sigma_2\rangle \cap \langle \sigma_3
		\rangle\;,\\
	\langle \sigma_1,\sigma_2 \rangle \cap
	\langle \sigma_2,\sigma_3 \rangle& =& \langle \sigma_{2} \rangle .

   \end{array}
\end{equation}

Conversely, if a group 
$\Lambda = \langle \sigma_1 , \sigma_2 , \sigma_3 \rangle$ 
satisfies (\ref{chirrels}) and  (\ref{chirinter}),
then there
exists a chiral or directly regular 4-polytope 
$\mathcal{P} = \mathcal{P}(\Lambda)$ 
 of type $\{ p_1,p_2,p_3\}$. 
 We refer to \cite{mowe} or \cite[Thm. 1]{schwei} 
 for the details of the construction.
The directly regular case occurs if and only if $\Lambda$ admits 
an involutory automorphism $\rho$
such that $(\sigma_1)\rho=\sigma_1^{-1}$, 
$(\sigma_2)\rho = \sigma_1^{2}\sigma_2$
and $(\sigma_3)\rho = \sigma_3$.

A chiral polytope $\mathcal{P}$ can be 
 self-dual in two subtly different ways
(see \cite{hub} or \cite[\S\,3]{scwe}).  
 $\mathcal{P}$ is \textit{properly} self-dual 
if it admits a
polarity $\delta$ which reverses the base flag $\Phi$ and
so preserves the two flag orbits. In $D(\mathcal{P})$ we  
then have 
$\delta^2 = \epsilon$ and $\delta \sigma_j \delta = \sigma_{4-j}^{-1}$, 
for $j = 1,2,3$.
In contrast,  $\mathcal{P}$ is {\em improperly} self-dual if there exists a 
duality $\delta$  which
exchanges the two flag orbits. In fact, we can choose $\delta$ so that  
$\delta^2 = \sigma_1 \sigma_2 \sigma_3$ (so $\delta$ has period $4$); and
$\delta^{-1} \sigma_1 \delta= \sigma_3^{-1},\,
\delta^{-1} \sigma_2 \delta = \sigma_1 \sigma_2 \sigma_1^{-1},\,
\delta^{-1} \sigma_3 \delta = \sigma_1$.

In Theorem 5 of  \cite{mowe}
we find that $\mathcal{G}$ is 2-transitive when $\mathcal{P}$ is finite,
chiral and self-dual of type $\{3,q,3\}$; more specifically, $\mathcal{G}$
is then of type $2^+$ (resp. $2^-$)
if and only if $\mathcal{P}$ is properly (resp. improperly) self-dual.

\medskip

In the above results, the self-duality of the polytope $\mathcal{P}$
serves as a natural guarantee that the medial layer graph $\mathcal{G}$ 
be vertex-transitive.
Now ignoring duality, it is quite clear from the symmetry of 
$\mathcal{P}$ that $\rm{Aut}(\mathcal{G})$ is transitive on the edges
of $\mathcal{G}$, and separately, at least, is also transitive
on $t$-arcs of types $1$ or $2$, for some $t \geq 2$.
We thus ask whether $\rm{Aut}(\mathcal{G})$ can be transitive on
all such $t$-arcs, thereby making $\mathcal{G}$ 
\textit{symmetric}, even when
$\mathcal{P}$ is \textit{not} self-dual. In fact, we shall see that 
this cannot happen, and so we make the following

\begin{definition} A finite regular
 graph $\mathcal{G}$ is {\em semisymmetric}
if $\rm{Aut}(\mathcal{G})$ acts transitively on the edges of
$\mathcal{G}$ but not transitively on the vertices of $\mathcal{G}$.
\end{definition}

\noindent
\textbf{Remarks}. To be quite clear about terminology, we recall that
a `regular' graph has all 
vertices of some fixed degree $k$. 
Semisymmetric graphs are a little elusive and 
hence of considerable interest.  The so-called \textit{Gray graph}
is the earliest known example of a trivalent semisymmetric graph;
see \cite{biggs},  \cite{bou1,bou2},  or \cite{marpar,pisaram} 
for neat descriptions, and
\cite{condmal} for
another interesting `small' example. A census of such graphs, with at most
$768$ vertices, appears in \cite{condsemi}.

It is easy to check that  a connected, semisymmetric graph 
$\mathcal{G}$ is bipartite, say with vertices of types $1$ and $2$. 
In analogy to the symmetric case, we define $\mathcal{G}$
to be \textit{$(t_1, t_2)$-semitransitive} if, for $j=1,2$,
${\rm Aut}(\mathcal{G})$ is transitive on $t_j$-arcs emanating from vertices
of type $j$ (but of course not transitive on longer such arcs).
In brief, we   say then that $\mathcal{G}$ is 
\textit{ss of type $(t_1,t_2)$}. 
The  theory of such graphs seems to be largely
uncharted, although it was proved in \cite{rweiss1}
that each $t_j \leq 7$. A further generalization is the notion 
of a \textit{locally $s$-arc transitive graph}; see   
 \cite{giud1}  for a detailed survey, or  \cite{giud2,giud3}
for more specific investigations. We note that the
`$s$-arc transitivity'  discussed in the papers
just cited has a more general meaning than that employed here.

In the next section we develop some machinery for manufacturing
semisymmetric trivalent graphs from non-self-dual
regular or chiral  $4$-polytopes of type $\{3,q,3\}$.

\medskip 
 
\section{Vertex-transitive medial layer graphs}

We begin by letting   $\CP$ be a regular polytope of type $\{3,q,3\}$,
with  medial layer graph $\CG$.   Theorem~\ref{vertransreg} below
characterizes the case in which $\rm{Aut}(\CG)$
is vertex-transitive.

\begin{theorem}
\label{vertransreg}
Suppose that $\CP$ is a finite regular $4$-polytope of type
$\{3,q,3\}$ with medial layer graph $\CG$. Then if
$\CG$ is vertex-transitive, $\CG$ must actually be $3$-transitive and
$\CP$ must be self-dual.
\end{theorem}

\textbf{Proof}.
Suppose that $\CG$ is vertex-transitive. Then $\CG$ must be
transitive on $3$-arcs. In fact, $\rm{Aut}(\CG)$ is
already known to be transitive on the $3$-arcs of each type, and any
element of $\rm{Aut}(\CG)$ which maps a
vertex $x$ of $\CG$ to a vertex $s$ of different type must
necessarily also map a $3$-arc with initial vertex
$x$ to a $3$-arc of the other type, with initial vertex $s$.

Next we show that $\rm{Aut}(\CG)$ is actually sharply transitive on
$3$-arcs, that is, $\CG$ is $3$-transitive.
We need to exclude the possibility that $\CG$ is $t$-transitive for
$t=4$ or $5$. In the notation of the previous section, we now have
$\rm{Aut}(\mathcal{P}) = \langle\rho_{0},\ldots,\rho_{3}\rangle$. It is
also useful to specify a few more vertices in Figure~\ref{frag}:
let $x:=(v_{-1})\rho_3$ and $y$ to be the two 
other  vertices adjacent to $w_1$, and likewise let
$s:=(v_4)\rho_0$ and $z$ be the two others adjacent to $w_2$.

The case $t=4$ can be ruled out as in \cite[Thm. 2]{mowe}, using
the fact that the stabilizer $B_4$ of a
vertex in a finite connected $4$-transitive trivalent graph must be
isomorphic to $\mathbb{S}_4$. In fact, the
element $\eta :=\rho_{0}\rho_{2}\rho_{3}$ in $\rm{Aut}(\CP)$ is an
automorphism
of $\CG$ that stabilizes the
vertex $v_{1}=F_{1}$ of $\CG$ and has order $6$; and it permutes the
vertices at distance $2$ from $v_1$ in the
$6$-cycle $(x\,w_{0}\,v_{3}\,y\,v_{-1}\,w_{2})$. 
However, $\mathbb{S}_4$ does not contain an element of order $6$.

The elimination of the case $t=5$ is more elaborate. When $t=5$, the
stabilizer $B_5(v_1)$ of the vertex
$v_{1}=F_{1}$ of $\CG$ in $\rm{Aut}(\CG)$ must be isomorphic to $\mathbb{S}_4
\times \mathbb{Z}_2$. However, the
stabilizer of $v_1$ in $\rm{Aut}(\CP)$ is just the subgroup
$\langle\rho_{0},\rho_{2},\rho_{3}\rangle \cong
\mathbb{S}_{3} \times \mathbb{Z}_2$. We claim that $\rho_0$ is the central
element of $B_5(v_1)$, determining the
factor $\mathbb{Z}_2$. In fact, viewing $\mathbb{S}_4\times \mathbb{Z}_2$ as
the symmetry group $[4,3]$ of the
$3$-cube $\{4,3\}$, we observe that its only subgroups of type $\mathbb{S}_3$
are those that fix a vertex of the cube,
and that the central inversion is the only non-trivial element in
$[4,3]$ that commutes with a subgroup of
this kind. Hence $\rho_0$, which determines the factor $\mathbb{Z}_2$
in $\mathbb{S}_{3}\times \mathbb{Z}_2$, is the
central element of $B_5(v_1)$.

Since the vertex-stabilizers in $\rm{Aut}(\CG)$ are all conjugate, we
similarly obtain that $\rho_1$,
$\rho_2$ and $\rho_3$ are the central elements in the stabilizers of
the vertices $w_2$, $w_1$ and
$v_{2}=F_2$, respectively, denoted by $B_5(w_2)$, $B_5(w_1)$ and
$B_5(v_{2})$ (see Figure~\ref{frag}). Now consider an element 
$\delta$ in $\rm{Aut}(\CG)$ which
maps the $3$-arc
$[w_{1},v_{1},v_{2},w_{2}]$ to the reversed $3$-arc
$[w_{2},v_{2},v_{1},w_{1}]$. 
Since the $\rho_j$'s are distinguished as central elements of their
respective vertex-stabilizers,
we   must therefore have
$\delta^{-1}\rho_j\delta = \rho_{3-j}$ for $j=0,1,2,3$.
Suppose for a moment that $\delta$ is an
involution. Then it follows that conjugation by $\delta$ in
$\rm{Aut}(\CG)$ induces an involutory group
automorphism of $\rm{Aut}(\CP)$, so that necessarily $\CP$ is
self-dual (see \cite[2E12]{arp}), contrary to our
assumption that $t=5$. (Recall from \cite[Thm. 2]{mowe}
that the medial layer graph of $\mathcal{P}$
must be
$3$-transitive if $\mathcal{P}$ is self-dual.)

It remains to prove that we may take $\delta$ to be an involution.
First observe that $\delta^2$ belongs
to the pointwise stabilizer of the $3$-arc
$[w_{1},v_{1},v_{2},w_{2}]$, which is isomorphic to
$\mathbb{Z}_2\times\mathbb{Z}_2$ (Theorem~\ref{autgraph}(b)).
In all cases,  $\delta^4 = \epsilon$ and $x\delta^2 = x$ or $y$. If
$x\delta^2 = x$, then $(x\delta)\delta^2 = x\delta$, so that 
$\delta^2$ fixes the $5$-arc $[x, w_1, v_1, v_2, w_2, x\delta ]$. Then
$\delta^2 = \epsilon $ as desired. Otherwise, $x\delta^2 = y$ and $y\delta^2 = x$.
Consider the unique
automorphism $\gamma$ of $\CG$ which fixes the $4$-arc
$[x,w_{1},v_{1},v_{2},w_{2}]$ pointwise but
interchanges $s$ and $z$ (i.e. $x\delta$ and $y\delta$). 
Then $\gamma\delta$ reverses the $3$-arc
$[w_{1},v_{1},v_{2},w_{2}]$, so we may
replace $\delta$ by $\gamma\delta$. But
$$ x(\gamma\delta)^2 = x\delta\gamma\delta  
= y\delta^2 = x , $$
so now $\gamma\delta$ is the desired involution.

It follows that $\CG$ is $3$-transitive. Now we apply the methods of
\cite[\S 4]{mowe}. In particular,
associated with $\CG$ is a certain subgroup $\Gamma$ of $\rm{Aut}(\CG)$
with a canonically defined
set of four involutory generators (see \cite[Def. 3]{mowe}), and
this subgroup $\Gamma$ is the automorphism
group of a certain self-dual ranked partially ordered set (see
\cite[Thm. 3]{mowe}). In the present context
we can actually identify the generators of $\Gamma$ with the
generators $\rho_j$ for $\rm{Aut}(\CP)$ (and hence
$\Gamma$ with $\rm{Aut}(\CP)$), and then also the new partially ordered set
with $\CP$ itself. Thus $\CP$ is
self-dual. This completes the proof.
\hfill $\Box$
\medskip

\medskip

The situation for chiral polytopes is quite similar. 
We give fewer details in the proof, which relies more closely 
on ideas used in establishing
\cite[Thm. 5]{mowe}.

\begin{theorem}
\label{vertranschir}
Suppose that $\CP$ is a finite chiral $4$-polytope of type
$\{3,q,3\}$ with medial layer graph $\CG$. Then if
$\CG$ is vertex-transitive, $\CG$ must actually be $2$-transitive and
$\CP$ must be self-dual.
\end{theorem}

\textbf{Proof}. 
Let $\CG$ be vertex-transitive. First observe that $\CG$ is
transitive on $2$-arcs. In fact, $\rm{Aut}(\CP)$
(and hence $\rm{Aut}(\CG)$) is transitive on the $2$-arcs of each type,
and the vertex-transitivity allows us
again to swap the two kinds of $2$-arcs. It follows that $\CG$ is
$t$-transitive for $t=2$, $3$, $4$ or
$5$. We must establish that $t=2$. We now have
$\rm{Aut}(\mathcal{P}) = \langle\sigma_{1},\sigma_{2},\sigma_{3}\rangle$.

Suppose first that $t=3$. We   apply the methods of \cite[\S 4]{mowe}
to prove that $\CP$ must actually be
regular, not chiral. In fact, because $\CG$ is $3$-transitive, we
again have a subgroup $\Gamma$ of
$\rm{Aut}(\CG)$ with canonically defined generators
$\rho_0,\rho_1,\rho_2,\rho_3$ (see \cite[Def. 3]{mowe}).
Consulting \cite[Lemma 1]{mowe} and its proof we find that the products
$\rho_0\rho_1,\rho_1\rho_2,\rho_2\rho_3$ can be identified with the
generators $\sigma_1,\sigma_2,\sigma_3$
of $\Gamma(\CP)$ acting on $\CG$. (It is crucial here that 
$\mathcal{G}$ be $3$-transitive.)
Moreover, the self-dual ``regular"
ranked poset (with a
flag-transitive action) associated with $\Gamma$ as in \cite[Thm.
3]{mowe} is actually isomorphic to $\CP$.
In fact, this poset can be defined completely in terms of the generators
$\rho_0\rho_1,\rho_1\rho_2,\rho_2\rho_3$ of the ``rotation subgroup"
$\Gamma^+$ of $\Gamma$ (see
\cite[p.510]{schwei}), that is, in terms of the generators
$\sigma_1,\sigma_2,\sigma_3$ of $\rm{Aut}(\CP)$.
However, the poset associated with $\sigma_1,\sigma_2,\sigma_3$ is
just $\CP$ itself. Hence $\CP$ must
be regular. It follows that we cannot have $t=3$.

To rule out the cases $t=4, 5$ we mimic part of the proof of 
\cite[Thm. 5]{mowe}, which utilized 
certain  universal relations  satisfied by  generators
of $\rm{Aut}(\mathcal{G})$, as described in \cite[\S\,1]{cond1}. 
In each case it is impossible to achieve 
($\sigma_2 \sigma_3)^2 = \epsilon$, given the other relations in
(\ref{chirrels}). (We note that Theorem 5 of \cite{mowe} has
almost the same hypotheses as here, except that
$\mathcal{P}$ is there assumed to be self-dual; but this self-duality is 
used only to guarantee 
that the medial layer graph $\mathcal{G}$ be vertex-transitive.)

Thus we must have $t=2$. Now, as in the proof of \cite[Thm. 6]{mowe},
the sharp transitivity of $\rm{Aut}(\mathcal{G})$ on $2$-arcs
enables a definition of a duality on $\mathcal{P}$, whether
$\mathcal{G}$ is of type $2^+$ or $2^-$\,: 
see \cite[eqns. (7) and (13)]{mowe}.
\hfill$\square$

\section{Graphs from polytopes of type $\{3,q,3\}$}

There is a wealth of finite trivalent semisymmetric graphs that are medial
layer graphs of regular or chiral polytopes $\CP$ of type $\{3,6,3\}$. 
Necessarily, by
Theorems~\ref{vertransreg} and \ref{vertranschir}, $\CP$ must not be self-dual.
However, before exploring such polytopes we must first review some key 
constructions.   

For any pair ${\bf s} = (s,t)$ of integers satisfying 
$s^2 + st + t^2 > 1$, the \textit{toroidal map}  
$\{3,6\}_{\bf s}$  has the structure of a finite
$3$-polytope (or polyhedron), usually chiral, but regular
just when $s t (s-t) = 0$. Referring to \cite[1D]{arp}, we merely note here
that $\{3,6\}_{\bf s}$ is obtained from the regular triangular
tessellation
$\{3,6\}$ of the Euclidean plane by
factoring out a suitable subgroup of the group of translation symmetries.
Taking $v = s^2 + st + t^2$, we find that $\{3,6\}_{\bf s}$
has $v$ vertices, $3v$ edges, $2v$ triangular facets and a rotation 
group $\langle \sigma_1 , \sigma_2 \rangle$ of order $6v$.
The toroidal map $\{6,3\}_{\bf s}$ can be constructed similarly and is dual to 
$\{3,6\}_{\bf s}$, both as a map on a compact surface and as an abstract 
polyhedron.

In any regular (or chiral) $n$-polytope $\mathcal{P}$,
all facets are isomorphic to a particular  ($n-1$)-polytope, say  $\mathcal{M}$; 
likewise each \textit{vertex-figure}  $\mathcal{N}$ (maximal section over a vertex in 
$\mathcal{P}$) is isomorphic to one  ($n-1$)-polytope $\mathcal{N}$. 
Conversely, given regular ($n-1$)-polytopes $\mathcal{M}$, $\mathcal{N}$, 
there may or may not 
exist a regular $n$-polytope $\mathcal{P}$ with facets $\mathcal{M}$
and vertex-figures $\mathcal{N}$; but if one such polytope exists, then there is 
a \textit{universal polytope} of this type,  denoted
$$ \{ \mathcal{M}\, ,\, \mathcal{N} \}\; ,$$
and from which all others  are obtained by identifications
 \cite[4A]{arp}. Somewhat more intricate results like this hold for
chiral polytopes \cite{schwei,scwe}.
\smallskip

\subsection{Medial layer graphs of finite universal polytopes.}\label{finuniv}
\smallskip 

Rephrasing the introductory remarks above, we observe that 
every (finite) regular polytope $\CP$ of type $\{3,6,3\}$ has certain
facets $\{3,6\}_{\bf s}$ and vertex-figures $\{6,3\}_{\bf t}$, with ${\bf
s} = (s^k,0^{2-k})$, ${\bf t} = (t^l,0^{2-l})$;  here, $s \geq 2$ if $k =
1$ and $s \geq 1$ if $k = 2$; likewise 
$t \geq 2$ if $l = 1$ and $t \geq 1$ if 
$l = 2$. In particular, $\CP$ is a quotient of the (generally infinite)
universal regular $4$-polytope
$$ \CP_{\bf s,t} := \{\{3,6\}_{\bf s},\{6,3\}_{\bf t}\} \;. $$
(See \cite[Section 11E]{arp} for details. 
In some cases, the only available 
construction for $\CP_{\bf s,t}$ is via the corresponding string C-group,
which in turn is naturally defined by a presentation encoding the local structure
of the polytope. We can expect no simple 
expression for the order of the group.)

For certain small parameter values these
universal polytopes are known to be finite; however, the finite polytopes
$\CP_{\bf s,t}$ have not yet been completely enumerated. Clearly,  if ${\bf
s} \neq {\bf t}$, then $\CP$ cannot be self-dual and hence its medial layer 
graph
is semisymmetric.

We list in Table~\ref{graphst} data for
the medial layer graphs $\CG_{\bf s,t}$ of those
universal polytopes $\CP_{\bf s,t}$ which are known to be finite; in the
last column we use `ss-$(t_1, t_2)$' or `$3^+$', respectively, 
 to indicate that $\CG_{\bf s,t}$ is semisymmetric of type $(t_1, t_2)$
or $3$-transitive. (The type of the last semisymmetric graph 
with  $40320$  vertices seems to be beyond brute force calculation
in GAP  \cite{gap}, for example.)
Recall that $N$ is the number of vertices.

\begin{table}[hbt]
\begin{center}
\begin{tabular}{||c|c|r|c||}
\hline
$\bf s$  & $\bf t$  &  $N\;\;\;\;$     & \rm{Transitivity type} \\     
\hline \hline
$(1,1)$ & $(1,1)$ & $18$    & $3^{+}$ \\  
\hline
$(1,1)$ & $(3,0)$ & $54$    & ss-(4,3) \\          
\hline
$(2,0)$ & $(2,0)$ & $40$    & $3^{+}$\\    
\hline
$(2,0)$ & $(2,2)$ & $120$   & ss-(3,3) \\       
\hline
$(3,0)$ & $(3,0)$ & $486$   & $3^{+}$ \\   
\hline
$(3,0)$ & $(2,2)$ & $6912$ & ss-(3,3) \\        
\hline
$(3,0)$ & $(4,0)$ & $40320$ &ss-(\rm{?,?}) \\       
\hline
\end{tabular}
\end{center}
\caption{The medial layer graphs of the known finite polytopes $\CP_{{\bf s},{\bf
t}}$.}
\label{graphst}
\end{table}

\smallskip
When ${\bf s} = {\bf t}$, the universal polytope $\CP_{\bf s,t}$ is  
self-dual and generally has many self-dual quotients. For example,  
when the standard representations of the crystallographic Coxeter  
groups $[3,6,3]$ and $[3,\infty,3]$ are reduced modulo an odd prime $p 
$, we obtain interesting self-dual (in one case, non-self-dual)  
regular polytopes of types $\{3,6,3\}$ or $\{3,p,3\}$, respectively,  
with automorphism groups isomorphic to finite reflection groups over  
the finite field $\mathbb{Z}_p$ (see \cite[(28),(31)]{mosc}). The  
exception occurs for $[3,6,3]$ with $p=3$, yielding the non-self-dual  
polytope $\CP_{(1,1),(3,0)}$, whose medial layer graph is the Gray graph  
(see Section~\ref{graypoly}). All other polytopes 
obtained by this construction 
have finite trivalent  
symmetric graphs as medial layer graphs. In particular, when $p > 3$, the  
polytopes obtained from $[3,6,3]$ have facets $\{3,6\}_{(p,0)}$ and  
vertex-figures $\{6,3\}_{(p,0)}$ and hence are quotients of  $\CP_{(p, 
0),(p,0)}$.

\smallskip

\subsection{Non-constructive methods.}\label{nonconve}
\smallskip

Even if the universal polytope $\CP_{{\bf s},{\bf t}}$ 
of \S\ref{finuniv} is not finite, we
often can still establish the existence of semisymmetric medial layer graphs
through non-constructive methods by appealing to \cite[Thm. 4C4]{arp}.
Recall that a group $\Gamma$ is {\em residually finite} if, for each finite
subset of
$\Gamma\setminus\{\epsilon\}$, there exists a homomorphism of $\Gamma$ onto
a finite group such that no element of the subset is mapped to the identity
element.

Suppose that $Q$ is an infinite regular $4$-polytope with facets
$\{3,6\}_{\bf s}$ and vertex-figures
$\{6,3\}_{\bf t}$, whose group $\Gamma(\CQ)$ is residually finite. Then
\cite[Thm. 4C4]{arp}, applied
with $\CP_{1} = \{3,6\}_{\bf s}$ and $\CP_{2} = \{6,3\}_{\bf t}$, says that
there are infinitely many finite regular $4$-polytopes with facets
$\{3,6\}_{\bf s}$ and vertex-figures $\{6,3\}_{\bf t}$, which are quotients
of $\CQ$. When ${\bf s} \neq {\bf t}$, these polytopes yield trivalent
semisymmetric graphs.

Such polytopes $\CQ$ are known to exist at least for certain parameter
values, including
${\bf s} = (s,s)$ and ${\bf t} = (s,0)$ or $(3s,0)$, with $s\geq 2$ (but
excluding the pair ${\bf s}=(2,2)$ and
${\bf t}=(2,0)$). In fact, inspection of the methods employed in the proof
of \cite[Thm. 11E5]{arp} reveals the existence of certain infinite regular
$4$-polytopes $\CQ$ with facets $\{3,6\}_{\bf s}$ and vertex-figures
$\{6,3\}_{\bf t}$, whose group $\Gamma(\CQ)$ is a semi-direct product of an
infinite, finitely generated, $4$-dimensional complex linear group by a
small group ($\mathbb{S}_3$, in fact); 
then $\Gamma(\CQ)$ itself also is a complex linear
group, in a space of dimension larger than $4$ (see \cite[pp.
415-416]{arp}). By a theorem of Malcev~\cite{mal}, every finitely generated
linear group is residually finite. Thus $\Gamma(\CQ)$ is residually finite.

In summary, we obtain the following

\begin{theorem}
\label{noncon}
Let $s\geq 2$, and let ${\bf s} := (s,s)$ and ${\bf t} := (s,0)$ or
$(3s,0)$, but excluding the pair ${\bf s}=(2,2)$ and
${\bf t}=(2,0)$. Then there are infinitely many finite trivalent
semisymmetric graphs which are medial layer graphs of finite regular polytopes
with facets $\{3,6\}_{\bf s}$ and vertex-figures $\{6,3\}_{\bf t}$.
\end{theorem}

\medskip

As a final application of  these methods, 
we mention a similar such theorem for symmetric graphs.

\begin{theorem}
\label{nonconsym}
For each $q\geq 5$, there are infinitely many finite, trivalent  
symmetric (indeed  $3$-transitive)
graphs which are medial layer graphs of finite self-dual regular  
polytopes of type $\{3,q,3\}$.
\end{theorem}

\textbf{Proof}.
The Coxeter group $[3,q,3]$ is the automorphism group of the self- 
dual universal regular polytope $\CP:=\{3,q,3\}$.  In particular, $D 
(\CP) \cong [3,q,3] \rtimes \mathbb{Z}_2$, where $\mathbb{Z}_2$ is  
generated by the polarity $\delta$ that fixes the base flag of $\CP$  
($\delta$ corresponds to the symmetry of the string Coxeter diagram).  
Hence $D(\CP)$ is residually finite, since $[3,q,3]$ is residually  
finite. Now adapt the proofs of \cite[Thm. 4C4]{arp} and \cite[Cor.  
4C5]{arp}, applying Malcev's theorem to $D(\CP)$ in place of  $\rm 
{Aut}(\CP)$, and requiring that $\delta$ does not become trivial  
under the homomorphisms onto finite groups. Then the latter  
guarantees the self-duality of the resulting quotients of $\CP$;  
hence their medial layer graphs are symmetric.
\hfill $\Box$

\smallskip

\subsection{Polytopes and graphs from the Eisenstein integers.}\label{eisenpoly}
\smallskip 

Next we consider from \cite[\S\,6]{eis} a family of  regular or chiral
polytopes
$\mathcal{Q}_m^A$, again of type $\{3,6,3\}$.
Here, the parameter $m$ 
is chosen from $\mathbb{D} := \mathbb{Z}[\omega]$, the domain
of Eisenstein integers. (Recall that  $\omega = e^{2\pi i /3}$ is a primitive 
cube root of unity.) The construction begins with a certain 
group $H_m$ of $2 \times 2$ matrices over the residue class ring
$\mathbb{D}_m := \mathbb{D}/(m)$; and then any subgroup $A$   of the unit
group of $\mathbb{D}_m$, with $ -1 \in A$, is
said to be
\textit{admissible}. Without going into
many details, we note simply that the rotation group
$H_m^A = \langle \sigma_1, \sigma_2, \sigma_3 \rangle$
for $\mathcal{Q}_m^A$ is obtained from the matrix group by factoring 
out the subgroup consisting of scalar multiples of the identity, with
scalars from $A$. Thus $\mathcal{Q}_m^A$ is finite when $m \neq 0$.
On the other hand, $H_0^{\pm 1}$ is the infinite rotation group
for $\mathcal{Q}_0^{\pm 1}$, which is isomorphic to 
the regular honeycomb $\{3,6,3\}$ 
of hyperbolic space $\mathbb{H}^3$.

If the Eisenstein prime $1-\omega$ does not divide $m$, the polytope will be 
self-dual. Interesting as it is, we leave this case 
behind (see \cite{mowe}). Suppose therefore that
$$m = (1-\omega)^e d\;, $$
where $e \geq 1$ and $d \in \mathbb{D}\setminus\{0\}$. To avoid degeneracies,
we also assume that $d$ is a non-unit if $e=1$.
It follows from \cite[Thm. 6.1]{eis} that $\mathcal{Q}_m^A$ is 
a finite quotient of the universal polytope
$$\{\, \{3,6\}_{(c,b)}\, , \, 
\{ 6,3\}_{(\frac{c-b}{3} , \frac{c+2b}{3})}\,\}\;, 
$$
where $m = c-b\omega$, for certain $b,c \in \mathbb{Z}$. 
Since the facets of $\mathcal{Q}_m^A$ 
are clearly not dual to its vertex-figures,   
$\mathcal{Q}_m^A$ itself cannot be self-dual. Furthermore, 
$\mathcal{Q}_m^A$ is regular if  
$m \mid \bar{m}$ and $A = \bar{A}$ (i.e. the scalar subgroup
is  invariant under complex conjugation); 
otherwise, $\mathcal{Q}_m^A$
is chiral.
Consequently, by Theorems ~\ref{vertransreg} and \ref{vertranschir} above,
we obtain a trivalent, semisymmetric medial layer 
graph $\mathcal{G}_m^A$ with
$$
N = 2 \,[\,\frac{(m \bar{m})^3}{12\cdot |A|}\; \prod_{\pi|m} (1 -
(\pi \bar{\pi})^{-2})\,]
$$
vertices. (The product here is over all \textit{non-associated} prime 
divisors $\pi$ of $m$.) We can summarize our construction in the following 

\begin{theorem}
Suppose the Eisenstein integer $m$   satisfies $m\bar{m} = 3k$, for some
rational integer $k > 1$; and let $A$ be any admissible group of scalars. Then
$\mathcal{G}_m^A$ is a finite, trivalent semisymmetric graph.
\end{theorem}

\noindent
\textbf{Remarks}. When $m = 3 = (1-\omega)^2(-\omega^2)$
and $A = \{ \pm 1\}$, we get  the dual of the universal polytope
$\{ \{3,6\}_{(1,1)} , \{6,3\}_{(3,0)} \}$ mentioned in
\S \ref{finuniv} above. The medial layer graph $\mathcal{G}_3^{\pm 1}$
is the Gray graph, which we examine more closely below. (For easier  reading 
we  omit  the  brackets from $\{\pm 1\}$.)  Similarly, for $m = 2(1-\omega)$
we find that $\mathcal{Q}_{2-2\omega}^{\pm 1}$ is the dual
of the universal polytope 
$\{ \{3,6\}_{(2,0)} , \{6,3\}_{(2,2)} \}$ described in \S \ref{finuniv}.
The medial layer graph has $120$ vertices. According
to the census in \cite{condsemi},  we have thus
described the unique trivalent semisymmetric graphs with these
orders.

We note that $\mathcal{Q}_m^A$ itself is 
not usually the universal polytope for the specified toroidal facets and
vertex-figures.
Certainly we get a proper quotient of the universal
cover when $|A| > 2$, which is possible
when $m$  has distinct prime divisors. In any case, the scalar group $A$  
has order $2^a$ and depends in 
an intricate way on the prime factorization of $m$
in $\mathbb{D}$; see \cite[pp. 105-106]{eis}.

\subsection{The Gray graph.}\label{graypoly}
\smallskip 

The Gray graph $\mathcal{C}$ is the smallest  
trivalent, semisymmetric graph (see \cite{condsemi}). Following  
\cite{bou2}, we define
$\mathcal{C}$ to be the (bipartite) incidence graph of cubelets and columns
in  a $3 \times 3 \times 3$ cube. Thus  vertices of the first type 
are the $27$ cubelets; and  vertices of the second type are the $9+9+9$
columns of $3$ cubelets parallel to edges of the cube. It is not hard to check
that $|\rm{Aut}(\mathcal{C}| = 1296$ \cite[Thm. 1.1]{bou2}.
Recent work has concerned
various interesting features of the graph
(\cite{marpar} and \cite{marusic-pisanski-wilson}); and here, of course, 
we construct it in a new way.
  
Before confirming that $\mathcal{G}_3^{\pm 1}$ 
\textit{really is} isomorphic to $\mathcal{C}$, 
we    develop a more concrete geometric 
description. First of all, using  
GAP  it is easy to check that 
$|\rm{Aut}(\mathcal{G}_3^{\pm 1})| = 1296$. Somewhat unexpectedly
this is $4$ times the order of
$ \rm{Aut}(\mathcal{Q}_3^{\pm 1})$.  
This discrepancy hints that we might 
examine a related embedding of the
honeycomb $\{3,6,3\}$ into a different  hyperbolic honeycomb $\{3,3,6\}$.  

\smallskip

\begin{figure}[hbt]
\begin{center}
\includegraphics[width=90mm]{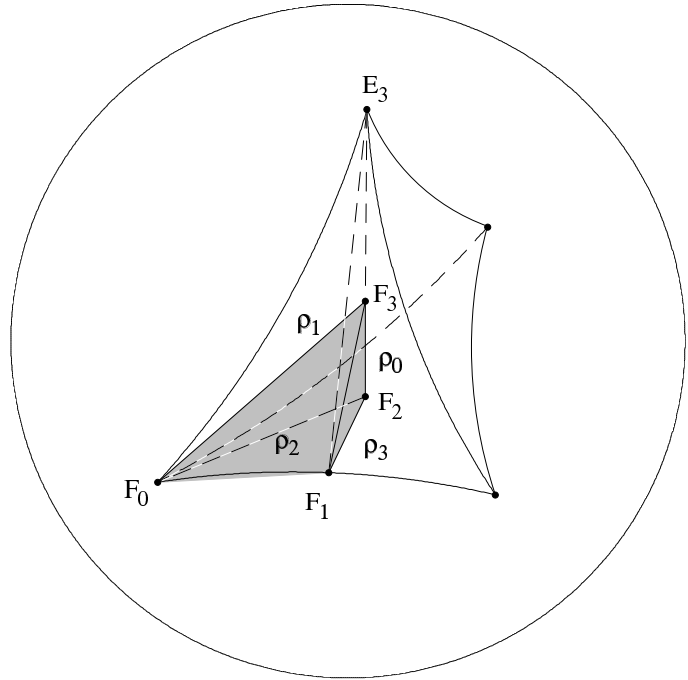}
\end{center}
\caption{A tetrahedral tile in the hyperbolic honeycomb $\{3, 3, 6 \}$.}
\label{Ttile}
\end{figure}

\smallskip
\newpage

Let us pick an arbitrary tetrahedral tile $T$ of  $\{3, 3, 6\}$ 
and denote its 
centre by $F_3$. In Figure~\ref{Ttile}, $F_0$ is a vertex of $T$
(and is an ideal point on the sphere at infinity);   let $F_1$ be the centre of 
an edge of $T$ through $F_0$  
and $ F_2$ the centre of a triangle of $T$ with that edge. The points 
$F_i$ ($i=0, 1, 2, 3$) are the vertices of a 
fundamental region for the hyperbolic Coxeter group
$[3, 3, 6]$. Thus, taking $\rho_i$ to be  the reflection in the face  
opposite  $F_i$ in this fundamental region, we have
$[3,3,6] = \langle \rho_0, \rho_1, \rho_2, \rho_3\rangle$.

Now let $E_3$ be the vertex of the tetrahedron $T$ 
which does not 
belong to the triangle centred at $F_2$. Then $F_0, F_1, F_2$ and $E_3$ 
are vertices of a new fundamental 
region for the hyperbolic Coxeter group
$[3, 6, 3]$, appearing here as the subgroup of $[3,3,6]$ generated by
$$ \rho_0, \rho_1, \rho_2^{\prime} := \rho_2\rho_3\rho_2\; \rm{ and }\; \rho_3\;.
$$
In fact, $[3,6,3]$ has index $4$ in $[3,3,6]$ (see \cite[11G]{arp}).

\begin{figure}[hbt]
\begin{center}
\includegraphics[width=60mm]{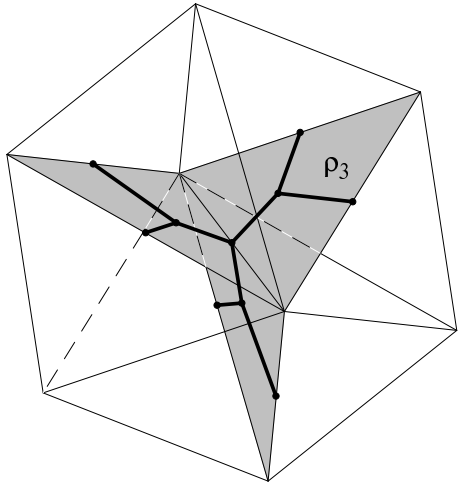}
\end{center}
\caption{The medial layer graph of $\{3, 6, 3\}$ inscribed in $\{3, 3, 6 \}$.}
\label{medinscr}
\end{figure}

Through the edge containing $F_1$ there are six triangles of $\{3, 3, 6\}$,
but only three of them 
belong to the honeycomb $\{3, 6, 3\}$, as we indicate  in Figure~\ref{medinscr}. 
Hence, the vertices of the medial layer graph 
of $\{3, 6, 3\}$ are comprised of just `half' the edges and `one quarter' of 
the triangles of $\{3, 3, 6\}$. 

\medskip

\begin{figure}[hbt]
\begin{center}
\includegraphics[width=60mm]{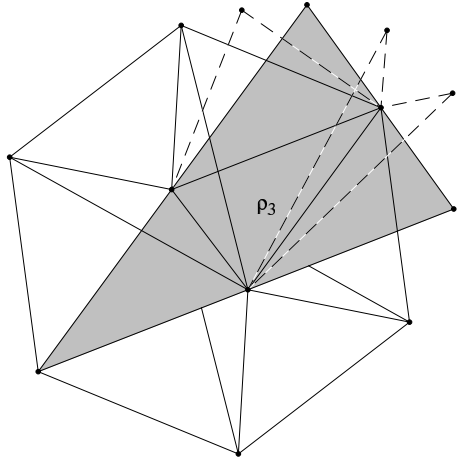}
\end{center}
\caption{The tessellation  $\{3,\infty\}$ cutting through  $\{3, 3, 6 \}$.}
\label{tesscut}
\end{figure}

\medskip 
Now a typical $2$-face of $\{3, 6, 3 \}$ is the ideal triangle $\{3\}$
with vertex $F_0$, edge $F_1$ and center $F_2$ in the hyperbolic
plane $p$ which serves as the mirror for the reflection $\rho_3$.
But $p$  is perpendicular to the mirrors for reflections
$\rho_0, \rho_1$ and 
$\rho:=\rho_2^{\prime} \rho_2 \rho_2^{\prime} = \rho_2(\rho_3\rho_2)^2$
in $[3,3,6]$. Since the latter two mirrors are parallel at $F_0$,
we see that $\langle \rho_0, \rho_1, \rho \rangle \simeq [3,\infty]$.
Hence the honeycomb $\{3, 3, 6 \}$ cuts the plane $p$ into the triangles of the
regular tessellation $\{3, \infty \}$ (see Figure~\ref{tesscut}).
 
\medskip

\medskip  
At this point we recall from \S\ref{eisenpoly} that the polytopes
$\mathcal{Q}_m^A$ ultimately arise from a modular representation
of the rotation subgroup in the reflection group $[3,6,3]$ for $\{3,6,3\}$.
A parallel construction based on the group $[3,3,6]$ is described
in \cite[\S\,4]{eis}, thereby yielding a family of polytopes 
$\mathcal{P}_m^A$ of type $\{\, \{3,3\}\, , \, \{3,6\}_{(b,c)}\,\}$.
When $1-\omega$ divides $m$, the rotation group $H_m^A$ for
 $\mathcal{Q}_m^A$ still has index $4$ in the rotation group 
$G_m^A$ for $\mathcal{P}_m^A$. In particular, when $m=3$ we obtain 
the universal polytope
$$\mathcal{P}_3^{\pm 1} = \{\, \{3,3\}\, , \, \{3,6\}_{(3,0)}\,\}\;,$$ 
whose automorphism group has order $1296$ \cite[11B5]{arp}. (Note that the 
rotation group for $\mathcal{Q}_3^{\pm 1}$ has order just $324$.) 
To see that the larger  group of order $1296$ is the  
automorphism group for the Gray graph, we must consider how reduction
modulo $m=3$  affects the picture in $\mathbb{H}^3$.
 
In identifying vertices of  $\{3, 3, 6\}$ to obtain the universal polytope 
$\mathcal{P}_3^{\pm 1}$, 
we note that $\rho_1 \rho = \rho_1 \rho_2 (\rho_3 \rho_2)^2$ has 
order $3$, so that the tessellation $\{3,\infty\}$ in the plane $p$
collapses to a `medial tetrahedron' $\{3,3\}$. 
At  the same time, the `inscribed' honeycomb 
$\{3, 6, 3\}$ collapses to $\mathcal{Q}_3^{\pm 1}$ and its medial layer 
graph to 
 $\mathcal{G}_3^{\pm 1}$.  

\begin{figure}[hbt]
\begin{center}
\includegraphics[width=60mm]{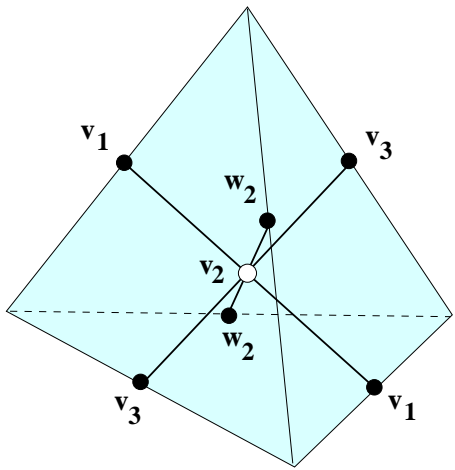}
\end{center}
\caption{The graph constructed from the `medial tetrahedra'.}
\label{medtetr}
\end{figure}

Now working in $\mathcal{P}_3^{\pm 1}$ , we define a graph $\mathcal{M}$ 
whose vertex set consists  of all  
$27 =  1296/(2\cdot 24)$ of the medial tetrahedra from
$\mathcal{P}_3^{\pm 1}$, together with all 
$ 27 =  (27\cdot 3)/3$  
pairs of opposite edges from
such tetrahedra. A vertex 
representing a tetrahedron is adjacent 
to a vertex representing a pair of edges whenever the tetrahedron contains 
the edges (see Figure~\ref{medtetr}). Evidently, $\mathcal{G}_3^{\pm 1}$ 
is isomorphic to $\mathcal{M}$ and so   
inherits all $1296$
automorphisms of $\mathcal{P}_3^{\pm 1}$.  In other words,

$$\rm{Aut}(\mathcal{G}_3^{\pm 1}) \cong \rm{Aut}(\mathcal{P}_3^{\pm 1})\;.$$
(This  group appears as $ [1 1 2]^3 \rtimes \mathbb{Z}_2$
in   \cite[11B5]{arp}.) 

\medskip

Now it is a simple matter to identify  
$\mathcal{G}_3^{\pm 1} \cong \mathcal{M}$ with 
the Gray graph $\mathcal{C}$: take the medial tetrahedra to
be the columns 
in the $3 \times 3 \times 3$ cube,
and pairs of opposite edges
to be  the  cubelets themselves. 
But each edge of $\{3,6,3\}$ belongs to three facets
$\{3,6\}$; so after reducing modulo $m = 3$, each of our $27$ pairs of edges 
must lie on three distint medial tetrahedra.

We note finally that it is not at all clear, for general 
moduli $m = (1-\omega)^e d$, just how large $\rm{Aut}(\mathcal{G}_m^A)$ is 
relative to its subgroup
$\rm{Aut}(\mathcal{Q}_m^A) $. 

For example, when $m = 2(1-\omega)$, 
$\mathcal{Q}_{2(1-\omega)}^{\pm 1}$ is a regular polytope whose
full reflection group of order $720$ is isomorphic to 
$\rm{Aut}(\mathcal{G}_{2(1-\omega)}^{\pm 1})$. A similar isomorphism holds
for the chiral polytope obtained when $m = (1-\omega)(1+3\omega)$. But when
$m = 3(1-\omega)$, the reflection group
$\rm{Aut}(\mathcal{Q}_{3(1-\omega)}^{\pm 1})$ once more has index $4$ in
the $\rm{Aut}(\mathcal{G}_{3(1-\omega)}^{\pm 1})$, whose order is $34992$. 

Based on this flimsy evidence, we conjecture that the index is always 
$4$ whenever
$m = (1-\omega)^e$, for $e\geq 2$.

\noindent
\textbf{Some history and words of thanks}. 
At this point we happily  thank Izak Bouwer for several
comments concerning the provenance of the Gray graph. 
In 1968, Izak gave the  first 
published description \cite{bou1}. A year later, 
in private 
correspondence with him, Dr. Marion C. Gray (1902--?)
wrote that
she had encountered the graph while
investigating `completely symmetric networks'. 
(This happened about 1932, early in her career at Bell Labs.)
In fact, 
with some uncertainty, 
Dr. Gray even attributed the graph to R. D. Carmichael.  
Perhaps the configurations described
in \cite{carm} were an inspiration.

It is also a pleasure to thank the referees for several
suggestions 
and for  pointing out related material in
\cite{giud1} and \cite{rweiss1}.



\end{document}